\documentclass{pnastwo}

\usepackage{amssymb,amsfonts,amsmath, tikz,caption}

\usepackage{color}

\newcommand{\cA}{\mathcal{A}}
\newcommand{\cL}{\mathcal{L}}
\newcommand{\cT}{\mathcal{T}}
\newcommand{\QQ}{\mathbb{Q}}
\newcommand{\ZZ}{\mathbb{Z}}

\newcommand{\uc}{\bar{e}}
\newcommand{\lc}{\underline{e}}
\newcommand{\UX}{\overline{X}}
\newcommand{\LX}{\underline{X}}

\contributor{Submitted to Proceedings
of the National Academy of Sciences of the United States of America}
\url{www.pnas.org/cgi/doi/10.1073/pnas.0709640104}
\copyrightyear{2008}
\issuedate{Issue Date}
\volume{Volume}
\issuenumber{Issue Number}

\begin{document}



\title{Greedy bases in rank 2 quantum cluster algebras}





\author{Kyungyong Lee\affil{1}{Department of Mathematics, Wayne State University, Detroit, MI 48202, U.S.A.},
            Li Li\affil	{2}{Department of Mathematics and Statistics, Oakland University, Rochester, MI 48309, U.S.A.},
            Dylan Rupel\affil{3}{Department of Mathematics, Northeastern University, Boston, MA 02115, U.S.A.}, \and
            Andrei Zelevinsky\affil{3}{}}

\contributor{Submitted to Proceedings of the National Academy of Sciences
of the United States of America}

\maketitle

\begin{article}

\begin{abstract} We identify a quantum lift of the greedy basis for rank 2 coefficient-free cluster algebras.  Our main result is that our construction does not depend on the choice of initial cluster, that it builds all cluster monomials, and that it produces bar-invariant elements.  We also present several conjectures related to this quantum greedy basis and the triangular basis of Berenstein and Zelevinsky \cite{triangular}. \end{abstract}








 \dropcap{C}luster algebras were introduced by Fomin and Zelevinsky \cite{fomin-zelevinsky1} as a combinatorial tool for understanding the canonical basis and positivity phenomena in the coordinate ring of an algebraic group.
 This was followed by the definition of quantum cluster algebras by Berenstein and Zelevinsky \cite{quantum} as a similar combinatorial tool for understanding canonical bases of quantum groups.

In connection with these foundational goals, it is important to understand the various bases of cluster algebras and quantum cluster algebras.  In this paper, we confirm the existence of a quantum greedy basis in rank 2 quantum cluster algebras. It is independent of the choice of an initial cluster,
contain all cluster monomials,
and specializes at $v=1$ to the greedy basis of a classical rank 2 cluster algebra introduced by Lee, Li, and Zelevinsky \cite{greedy}. 

Another important basis in an acyclic quantum cluster algebra is the triangular basis of Berenstein and Zelevinsky \cite{triangular}.  Like triangular bases, quantum greedy bases can be easily computed. This enables us to study and compare these bases computationally.

The structure of the paper is as follows. Section 1 begins with a recollection on rank 2 commutative cluster algebras.  We recall the construction of greedy bases here and in the process we provide a new axiomatic description of the greedy basis.  In Section 2 we review the definition of rank 2 quantum cluster algebras. In Section 3 we  present our main result: the existence of a quantum lift of the greedy basis using the axiomatic description.  Section 4 recalls the definition of the triangular basis and Section 5 provides an overview of results and open problems on the comparison between various bases of rank 2 quantum cluster algebras.



\section{1. Rank 2 Cluster Algebras and Their Greedy Bases}\label{sec:greedy}

Fix positive integers $b,c>0$.  The commutative cluster algebra $\cA(b,c)$ is the $\ZZ$-subalgebra of $\QQ(x_1,x_2)$ generated by the \emph{cluster variables} $\{x_m\}_{m\in\ZZ}$, where the $x_m$ are rational functions in $x_1$ and $x_2$ defined recursively by the \emph{exchange relations}
\[x_{m+1}x_{m-1}=\begin{cases} x_m^b+1 & \text{ if $m$ is odd;}\\ x_m^c+1 & \text{ if $m$ is even.}\end{cases}\]
It is a fundamental result of Fomin and Zelevinsky \cite{fomin-zelevinsky1} that, although the exchange relations appear to produce rational functions, one always obtains a Laurent polynomial whose denominator is simply a monomial in $x_1$ and $x_2$.  They actually showed the following slightly stronger result.

\begin{theorem}[{\cite[Theorem 3.1]{fomin-zelevinsky1}, Laurent Phenomenon}]
 For any\\\noindent $m\in\ZZ$ we have $\cA(b,c)\subset\ZZ[x_m^{\pm1},x_{m+1}^{\pm1}]$.
\end{theorem}

The cluster algebra $\cA(b,c)$ is of \emph{finite type} if the collection of all cluster variables is a finite set.  Fomin and Zelevinsky \cite{fomin-zelevinsky2} went on to classify cluster algebras of finite type.

\begin{theorem}[{\cite[Theorem 1.4]{fomin-zelevinsky2}}]
 The cluster algebra $\cA(b,c)$ is of finite type if and only if $bc\le3$.
\end{theorem}

\noindent The proof of this theorem in particular establishes a connection between denominator vectors of cluster variables and almost positive roots in a root system $\Phi\subset\ZZ^2$.  Thus we say $\cA(b,c)$ is of \emph{affine} (resp. \emph{wild}) \emph{type} if $bc=4$ (resp. $bc\ge5$).  

An element
$x\in\bigcap_{m\in\ZZ}\ZZ[x_m^{\pm1},x_{m+1}^{\pm1}]$ is called \emph{universally Laurent} since the expansion of $x$ in every cluster is Laurent.  If the coefficients in the Laurent expansion of a nonzero element $x$ in each cluster are positive integers, then $x$ is called \emph{universally positive}.  A universally positive element in $\cA(b,c)$ is said to be \emph{indecomposable} if it cannot be expressed as a sum of two universally positive elements.

Sherman and Zelevinsky \cite{sz-Finite-Affine} studied in great detail the collection of indecomposable universally positive elements.  One of their main results is the following: if the commutative cluster algebra $\cA(b,c)$ is of finite or affine type, then the indecomposable universally positive elements form a $\ZZ$-basis in $\cA(b,c)$, moreover this basis contains the set of all \emph{cluster monomials} $\big\{x_m^{a_1}x_{m+1}^{a_2}:a_1,a_2\in\ZZ_{\ge0},m\in\ZZ\big\}$.  

However the situation in the wild-types becomes much more complicated.  In particular, it was shown by Lee, Li, and Zelevinsky \cite{LLZ2} that for $bc\ge5$ the indecomposable universally positive elements of $\cA(b,c)$ are not linearly independent.  

The \emph{greedy basis} of $\cA(b,c)$ introduced by Lee, Li, and Zelevinsky \cite{greedy} is a subset of the indecomposable universally positive elements which admits a beautiful combinatorial description.  In this section we will describe a new axiomatic characterization of the greedy basis.

The elements of the greedy basis take on a particular form which is motivated by a well-known pattern in the initial cluster expansion of cluster monomials.  An element $x \in \mathcal{A}(b,c)$ is \emph{pointed} at $(a_1, a_2) \in \ZZ^2$ if it can be written in the form
\begin{equation*}
x=x_1^{-a_1}  x_2^{-a_2} \sum_{p,q \geq 0} e(p,q) x_1^{bp} x_2^{cq}	
\end{equation*}
with $e(0,0)=1$ and $e(p,q)=e(a_1,a_2,p,q) \in \ZZ$ for all $p,q\ge0$.  
It is well-known that $x_1^{a_1}  x_2^{a_2}$ is the denominator of a cluster monomial if and only if $(a_1,a_2)\in\ZZ^2\setminus \Phi^{im}_+$ (cf. \cite{CC,CK}), where
$\Phi^{im}_+$ is the set of positive imaginary roots, i.e.
\[\Phi^{im}_+:=\Big\{(a_1,a_2)\in\mathbb{Z}_{>0}^2: ca_1^2-bca_1a_2+ba_2^2\le 0\Big\}.\]

For $(a_1,a_2)\in\Phi^{im}_+$, define the region 
$$\aligned
R_{\text{greedy}}=&\bigg\{(p,q)\in\mathbb{R}_{\ge0}^2\bigg|\; q+\Big(b-\frac{ba_2}{ca_1}\Big)p<a_1 \textrm{\, or \,}\\
&p+\Big(c-\frac{ca_1}{ba_2}\Big)q<a_2\bigg\}
\cup\Big\{(0,a_1),(a_2,0)\Big\}.
\endaligned
$$
In other words, $R_{\text{greedy}}$ is the region bounded by the broken line

\hspace{30pt} $(0,0), (a_2,0), (a_1/b,a_2/c), (0,a_1), (0,0),$

\noindent with the convention that this region includes the closed segments $[(0, 0), (a_2 , 0)]$
and $[(0, a_1 ), (0, 0)]$ but excludes the rest of the boundary (see Fig. 1).

We can define greedy elements in two different ways, either by axioms or by recurrence relations. Here we choose to define them using recurrence relations as follows.

\begin{theorem}[{\cite[Proposition 1.6]{greedy}}] \label{recursive-definition}
For each $(a_{1},a_{2})\in\ZZ^{2}$, there exists a unique element in $\mathcal{A}(b,c)$ pointed at $(a_1,a_2)$ whose coefficients $e(p,q)$ satisfy the following recurrence relation:
$e(0,0)=1$,

$e(p,q)=$
$$  \begin{cases}
          \sum\limits_{k=1}^p(-1)^{k-1}e(p-k,q){ [a_2-cq]_++k-1\choose k} & \text{ if $ca_1q\le ba_2p$;}\\
          \sum\limits_{\ell=1}^q(-1)^{\ell-1}e(p,q-\ell){[a_1-bp]_++\ell-1\choose \ell} & \text{ if $ca_1q\ge ba_2p$,}
         \end{cases}
 $$where we use the standard notation $[a]_+=\max(a,0)$.
\end{theorem}

We define the {\it greedy element} pointed at $(a_1,a_2)$, denoted $x[a_{1},a_{2}]$, to be the unique element determined by Theorem \ref{recursive-definition}.

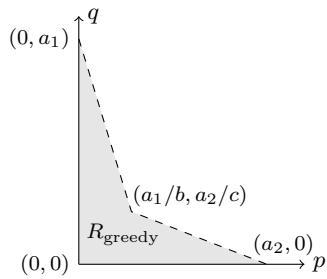
\begin{figure}
\begin{center}
\begin{tikzpicture}
\usetikzlibrary{patterns}
\fill [black!10]  (0,3)--(.7,.7)--(2.5,0)--(0,0)--(0,3);
\draw[dashed] (0,3)--(.7,.7)--(2.5,0);
\draw (0,0) node[anchor=east] {\small$(0,0)$};
\draw (2.2,0) node[anchor=south west] {\small$(a_2,0)$};
\draw (.6,.65) node[anchor=south west] {\small$(a_1/b,a_2/c)$};
\draw (0,3) node[anchor=east] {\small$(0,a_1)$};
\draw (1.2,0.4) node[anchor=east] {\small$R_{\text{greedy}}$};
\draw[->] (0,0) -- (3,0)
node[right] {$p$};
\draw[->] (0,0) -- (0,3.3)
node[right] {$q$};
\end{tikzpicture}
\end{center}
\caption{Support region of (quantum) greedy elements.}
\end{figure}

For  $(a_1,a_2)\in\Phi^{im}_+$, we give a new characterization for the greedy element.
\begin{theorem}\label{three axioms: commutative}
For each $(a_1,a_2)\in\Phi^{im}_+$, the pointed coefficients $e(p,q)$ of the greedy element $x[a_{1},a_{2}]$ are uniquely determined by the following two axioms:
\vspace{-.05in}
\begin{itemize} 
\item (Support) $e(p,q)=0$ for $(p,q)\not\in R_{\text{greedy}}$; \item (Divisibility)
\begin{itemize}
 \item if $a_2>cq$, then $(1+t)^{a_2-cq}\Big|\sum_i e(i,q)t^i$;
 \item if $a_1>bp$, then $(1+t)^{a_1-bp}\Big|\sum_i e(p,i)t^i$.
\end{itemize}\end{itemize}
\end{theorem}

The key step in the proof of Theorem~\ref{three axioms: commutative} is based on an observation made in \cite[Section 2]{greedy}: the first (resp. second) recurrence relation  in Theorem~\ref{recursive-definition} is equivalent to the vanishing of the $(p,q)$-th coefficient of the Laurent expansion of $x[a_1,a_2]$ with respect to the cluster $\{x_0,x_1\}$ (resp. $\{x_2,x_3\}$).  It is straightforward to check that such vanishing implies the Support and Divisibility axioms in Theorem~\ref{three axioms: commutative}.
\medskip

The following theorem summarizes the results from \cite{greedy}.

\begin{theorem}[{\cite[Theorem 1.7]{greedy}}]
\label{main theorem-commutative}\mbox{}

{\rm(a)} The greedy elements $x[a_1,a_2]$ for $(a_1, a_2) \in \ZZ^2$ form a $\ZZ$-basis in $\cA(b,c)$, which we refer to as the \emph{greedy basis}.

{\rm(b)} The greedy basis is independent of the choice of an initial cluster.

 {\rm({c})} The greedy basis contains all cluster monomials.

{\rm(d)}  Greedy elements are universally positive and indecomposable.
\end{theorem}

Our goal in this work is to generalize the above theorem to the setting of rank 2 quantum cluster algebras.  The proof of Theorem~\ref{main theorem-commutative} given in \cite{greedy} uses combinatorial objects called \emph{compatible pairs} in an essential way (cf. \cite[Theorem 1.11]{greedy}). Unfortunately this method has difficulties/limitations in generalizing to the study of quantum cluster algebras. More precisely, if one could simply assign a power of $v$ to each compatible pair then the quantum greedy elements would have again been universally positive, which is unfortunately false in general (see the example at the end of Section 3). This motivates our new approach to Theorem~\ref{main theorem-commutative}, with the exception of (d), 
which can easily be generalized to the quantum setting.

\section{2. Rank 2 Quantum Cluster Algebras}
 In this section we define our main objects of study, namely quantum cluster algebras, and recall important fundamental facts related to these algebras.  We restrict attention to rank 2 quantum cluster algebras where we can describe the setup in very concrete terms.  We follow (as much as possible) the notation and conventions of \cite{triangular,greedy}.

 Consider the following quantum torus
 $$\cT=\ZZ[v^{\pm 1}]\langle X_1^{\pm1},X_2^{\pm1}: X_2 X_1=v^2 X_1 X_2\rangle$$
 (this setup is related to the one in \cite{rupel1} which uses the formal variable $q$ instead of $v$ by setting $q = v^{-2}$).  There are many choices for quantizing cluster algebras, to rigidify the situation we require the quantum cluster algebra to be invariant under a certain involution.
The \emph{bar-involution} is the $\mathbb{Z}$-linear anti\-automorphism of $\cT$ determined by $\overline{f}(v)=f(v^{-1})$ for  $f\in\mathbb{Z}[v^{\pm1}]$ and
\[\overline{fX_1^{a_1}X_2^{a_2}}=\overline{f}X_2^{a_2}X_1^{a_1}=v^{2a_1a_2}\overline{f}X_1^{a_1}X_2^{a_2}\quad\text{($a_1,a_2\in\ZZ$)}.\]
An element which is invariant under the bar-involution is said to be \emph{bar-invariant}.

Let $\mathcal{F}$ be the skew-field of fractions
of $\cT$.    The \emph{quantum cluster algebra} $\cA_v(b,c)$ is the $\mathbb{Z}[v^{\pm1}]$-subalgebra of $\mathcal{F}$ generated by the \emph{quantum cluster variables} $\{X_m\}_{m \in\ZZ}$ defined recursively from the \emph{quantum exchange relations}
 \[X_{m+1}X_{m-1}=\begin{cases} v^{b}X_m^b+1 & \text{ if $m$ is odd;}\\ v^{c}X_m^c+1 & \text{ if $m$ is even.}\end{cases}\]

 By a simple induction one can easily check the following (quasi-)commutation relations between neighboring cluster variables $X_m,X_{m+1}$ in $\cA_v(b,c)$:
 \begin{equation}\label{eq:commutation-A11}
  X_{m+1} X_m = v^2 X_m X_{m+1} \quad (m \in \ZZ).
 \end{equation}
 It then follows that all cluster variables are bar-invariant, therefore $\cA_v(b,c)$ is also stable under the bar-involution.  Moreover, equation \eqref{eq:commutation-A11} implies that each \emph{cluster} $\{X_m, X_{m+1}\}$ generates a quantum torus
\[\cT_m=\ZZ[v^{\pm1}]\langle X_m^{\pm1},X_{m+1}^{\pm1}: X_{m+1} X_m = v^2 X_m X_{m+1}\rangle.\]
It is easy to see that the bar-involution does not depend on the choice of an initial quantum torus $\cT_m$.

 The appropriate quantum analogues of cluster monomials for $\cA_v(b,c)$ are the (bar-invariant) \emph{quantum cluster monomials} which are certain elements of  a quantum torus $\cT_m$, more precisely they are
 \[X^{(a_1,a_2)}_m = v^{a_1 a_2} X_m^{a_1} X_{m+1}^{a_2} \quad (a_1, a_2 \in \ZZ_{\geq 0}, \,\, m \in \ZZ).\]

\smallskip

The following quantum analogue of the Laurent Phenomenon was proven by Berenstein and Zelevinsky \cite{quantum}.

 \begin{theorem}[{\cite[Corollary 5.2]{quantum}}]\label{thm1}
  For any $m\in\ZZ$ we have $\cA_v(b,c)\subset\cT_m$.  Moreover, $$\cA_v(b,c)=\bigcap\limits_{m\in\ZZ}\cT_m.$$ 
 \end{theorem}

\smallskip

A nonzero element of $\cA_v(b,c)$ is called {\it universally positive} if it lies  in $\bigcap\nolimits_{m\in\ZZ}\ZZ_{\geq 0}[v^{\pm1}][X_m^{\pm1},X_{m+1}^{\pm1}]$. A universally positive element in $\cA_v(b,c)$ is {\it indecomposable} if it cannot be expressed as a sum of two universally positive elements.

\smallskip

 It is known as a very special case of the results of \cite{Q, Ef, DMSS, N-quiver, KQ} that  cluster monomials are universally positive Laurent polynomials.
Explicit combinatorial expressions for these positive coefficients can be obtained from the results of \cite{ls, rupel2, LL:qgrass}.
 \smallskip

\section{3. Quantum Greedy Bases}
This section contains the main result of the paper. Here we introduce the quantum greedy basis and present its nice properties.

Analogous to the construction of greedy elements, the elements of the quantum greedy basis take on the following particular form.
An element $X \in \mathcal{A}_{v}(b,c)$ is said to be {\it pointed} at $(a_1, a_2) \in \ZZ^2$ if it has the form
\begin{equation*}
X=\sum\limits_{p,q\ge0} e(p,q)X_1^{(bp-a_1,cq-a_2)}
\end{equation*}
with $e(0,0)=1$ and $e(p,q)\in\ZZ[v^{\pm1}]$ for all $p$ and $q$.
\medskip

The next theorem is an analogue of Theorem \ref{recursive-definition}, to state our result precisely we need more notation.  Let $w$ denote a formal invertible variable.  For $n\in\ZZ$, $k\in\ZZ_{\ge0}$ define the bar-invariant quantum numbers and bar-invariant quantum binomial coefficients by
 \[
 \aligned
 {}[n]_w&=\frac{w^n-w^{-n}}{w-w^{-1}}=w^{n-1}+w^{n-3}+\cdots+w^{-n+1};\\
 {n\brack k}_w&=\frac{[n]_w[n-1]_w\cdots[n-k+1]_w}{[k]_w[k-1]_w\cdots[1]_w}.
 \endaligned
 \]
 Recall that ${n\brack k}_w$ will always be a Laurent polynomial in $w$.  Hence taking $w=v^{b}$ or $v^{c}$ we obtain Laurent polynomials in $v$ as well.

\begin{theorem}\label{recursive-definition-q}
For each $(a_{1},a_{2})\in\ZZ^{2}$, there exists a unique element in $\mathcal{A}_v(b,c)$ pointed at $(a_1,a_2)$ whose coefficients $e(p,q)$ satisfy the following recurrence relation:
$e(0,0)=1$,

$e(p,q)=$
$$  \begin{cases}
          \sum\limits_{k=1}^p(-1)^{k-1}e(p-k,q){[a_2-cq]_++k-1\brack k}_{v^b} & \text{ if $ca_1q\le ba_2p$;}\\
          \sum\limits_{\ell=1}^q(-1)^{\ell-1}e(p,q-\ell){[a_1-bp]_++\ell-1\brack \ell}_{v^c} & \text{ if $ca_1q\ge ba_2p$.}
         \end{cases}
 $$
\end{theorem}
We define the {\it quantum greedy element} pointed at $(a_1,a_2)$, denoted $X[a_{1},a_{2}]$, to be the unique element determined by Theorem \ref{recursive-definition-q}. 

The following theorem is analogous to Theorem \ref{three axioms: commutative}. 

\begin{theorem}\label{three axioms: quantum}
For each $(a_1,a_2)\in\Phi^{im}_+$, the pointed coefficients $e(p,q)$ of $X[a_{1},a_{2}]$ are uniquely determined by the following two axioms:

\vspace{-.05in}
\begin{itemize}
\item (Support) $e(p,q)=0$ for $(p,q)\not\in R_{\text{greedy}}$;

\item (Divisibility) Let $t$ denote a formal variable which commutes with $v$.  Then
\begin{itemize}
 \item $\prod_{j=1}^{a_2-cq}\big(1+v^{b(a_2-cq+1-2j)}t\big)\Big|\sum_i e(i,q)t^i$;
 \item $\prod_{j=1}^{a_1-bp}\big(1+v^{c(a_1-bp+1-2j)}t\big)\Big|\sum_i e(p,i)t^i$.
\end{itemize}
\end{itemize}
\end{theorem}

Theorem \ref{three axioms: quantum} follows from an argument similar to the proof of Theorem \ref{three axioms: commutative}.

Our main theorem below states that quantum greedy elements possess all the desired properties described in Theorem \ref{main theorem-commutative} except the positivity (d).

\begin{theorem}
\label{main theorem}\mbox{}

{\rm(a)} The quantum greedy elements $X[a_1,a_2]$ for $(a_1, a_2) \in \ZZ^2$ form a $\ZZ[v^{\pm1}]$-basis in $\cA_v(b,c)$, which we refer to as the \emph{quantum greedy basis}.

{\rm(b)} The quantum greedy basis is bar-invariant and independent of the choice of an initial cluster.

 {\rm(c)} The quantum greedy basis contains all cluster monomials.

  {\rm(d)} If $X[a_1,a_2]$ is universally positive, then it is indecomposable.

  {\rm(e)} The quantum greedy basis specializes to the commutative greedy basis by the substitution $v=1$.
\end{theorem}
\smallskip

Full proofs of Theorems~\ref{recursive-definition-q}, \ref{three axioms: quantum} and \ref{main theorem} will appear elsewhere. 
The hardest part is to show the existence of quantum greedy elements, i.e. that the recursions in Theorem~\ref{recursive-definition-q} terminate. The main idea in the proof is to realize $X[a_{1},a_{2}]$ as a degeneration of certain elements in $\cA_v(b,c)$ whose supports are contained in the closure of ${R_{\text{greedy}}}$.

\medskip

We end this section with an example where positivity of quantum greedy elements fails.  Let $(b,c)=(2,3)$ and consider the pointed coefficient $e(2,1)$ in the Laurent expansion of $X[3,4]$.  Using Theorem~\ref{recursive-definition-q}, we see that $e(2,1)$ is not in $\mathbb{Z}_{\ge0}[v^{\pm1}]$:
$$\aligned
e(2,1)&=e(1,1){1\brack 1}_{v^2} -e(0,1){2\brack 2}_{v^2}\\
         &=e(1,1)-e(0,1)\\
         &=e(1,0){1\brack 1}_{v^3} -e(0,1)\\
         &=e(0,0){4\brack 1}_{v^2} - e(0,0){3\brack 1}_{v^3}\\
         &=(v^6+v^2+v^{-2}+v^{-6}) - (v^6+1+v^{-6})\\
         &=v^2-1+v^{-2}.
\endaligned$$
Further computation shows that  the greedy elements $X[3,5]$, $X[5,4]$, $X[5,7]$,  $X[5,8]$, $X[7,5]$, $X[7,10]$, $X[7,11]$, etc.~are not positive. Similarly,
positivity can be seen to fail for certain quantum greedy elements when $(b,c)=(2,5)$, $(3,4)$, $(4,6)$.

\section{4. Triangular Bases}\label{sec:triangular}


The construction of the triangular basis begins with the standard monomial basis.
 For every $(a_1,a_2) \in \ZZ^2$, we define the \emph{standard monomial} $M[a_1,a_2]$ (which is denoted $E_{(-a_1,-a_2)}$ in \cite{triangular}) by setting
 \begin{equation}\label{eq:Ma}
  M[a_1,a_2] = v^{a_1 a_2} X_3^{[a_1]_+} X_1^{[-a_1]_+} X_2^{[-a_2]_+} X_0^{[a_2]_+} \ .
 \end{equation}
 It is known \cite{quantum} 
 that the elements $M[a_1,a_2]$ form a $\ZZ[v^{\pm 1}]$-basis of $\cA_v(b,c)$.

 The importance of this basis comes from its computability, however this basis will not serve in our goals of understanding quantum cluster algebras.  Indeed, it is easy to see that the standard monomials are not bar-invariant and do not contain all the cluster monomials, moreover they are inherently dependent on the choice of an initial cluster. These drawbacks provide a motivation to consider the triangular basis
(as defined below) constructed from the standard monomial basis with a
built-in bar-invariance property.


 As the name suggests the triangular basis is defined by a triangularity property relating it to the standard monomial basis.  To describe this triangular relationship we introduce the ``$v$-positive'' lattice $\cL=\bigoplus v\ZZ[v]M[a_1,a_2]$ where the summation runs over all $a_1,a_2\in\ZZ$.  We  now define a new basis $\big\{C[a_1,a_2]\,:\,(a_1,a_2)\in\ZZ^2\big\}$ by specifying how it relates to the standard monomial basis:

 \begin{enumerate}
  \item[] (P1) Each $C[a_1,a_2]$ is bar-invariant.
  \item[] (P2) For each $(a_1,a_2) \in \ZZ^2$, we have
  $$\displaystyle{C[a_1,a_2] - M[a_1,a_2] \in \cL} \ .$$
 \end{enumerate}

 \smallskip

Using the universal acyclicity of rank 2 quantum cluster algebras, we  apply the following theorem.

\begin{theorem}[{\cite[Theorem 1.6]{triangular}}]
 The triangular basis does not depend on the choice of initial cluster and it contains all cluster monomials.
 \end{theorem}


Similar to the support condition for (quantum) greedy elements, we make the following conjecture about the supports of triangular elements.

\begin{conjecture}\label{conj:triangular support} Let $(a_1,a_2)\in\Phi^{im}_+$. For $0\le p\le a_2$, $0\le q\le a_1$, the pointed coefficient $e(p,q)$ of $X^{(bp-a_1,cq-a_2)}$ in $C[a_1,a_2]$ is nonzero if and only if 
\begin{equation*}
bp^2+bcpq+cq^2\le ca_1q+ba_2p.
\end{equation*}
\end{conjecture}

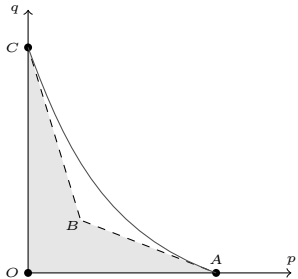
\begin{figure}[h]\label{fig:2}
\begin{center}
\begin{tikzpicture}
\usetikzlibrary{patterns}
\fill [black!10]  (0,3)--(.7,.7)--(2.5,0)--(0,0)--(0,3);
\draw[dashed] (.69,.73)--(0,3) (.73,.69)--(2.5,0);
\draw (0,0) node[anchor=east] {\tiny$O$};
\draw (2.5,0) node[anchor=south] {\tiny$A$};
\draw (.8,.8) node[anchor=north east] {\tiny$B$};
\draw (0,3) node[anchor=east] {\tiny$C$};
\draw[->] (0,0) -- (3.5,0)
node[above] {\tiny $p$};
\draw[->] (0,0) -- (0,3.5)
node[left] {\tiny $q$};
\fill (0,0) circle (1.5pt);
\fill (0,3) circle (1.5pt);
\fill (2.5,0) circle (1.5pt);
\draw [black!70] {(0,3) to [out=290, in=160] (2.5,0)};
\end{tikzpicture}
\end{center}
\caption{The conjectured support region of triangular basis elements is the closed region $OAC$ with a curved edge $AC$, where $O=(0,0)$, $A=(a_2,0)$, $C=(0,a_1)$. The support region of the (quantum or non-quantum) greedy element pointed at $(a_1,a_2)$ is the polygon $OABC$, where $B=(a_1/b,a_2/c)$. Note that
the line $\overline{BA}$ (resp.~$\overline{BC}$) is tangent to the curved edge $AC$ at point $A$ (resp.~$C$).
}
\end{figure}

\smallskip

 \section{5. Open Problems} \label{section5}

In this section we present open problems and conjectures on the triangular basis and quantum greedy basis  of a rank 2 quantum cluster algebra. 
 Our  aim is to find $\mathbb{Z}[v^{\pm1}]$-bases of $\cA_v(b,c)$ satisfying the properties laid out in the following definition.

A basis $\mathcal{B}$ for $\cA_v(b,c)$ is said to be {\it strongly positive} if the following hold:
\begin{enumerate}
\item[\rm(1)] each element of $\mathcal{B}$ is bar-invariant;
\item[\rm(2)] $\mathcal{B}$ is independent of the choice of an initial cluster;
\item[\rm(3)] $\mathcal{B}$ contains all cluster monomials;
\item[\rm(4)] any product of elements from $\mathcal{B}$ can be expanded as a linear combination of elements of $\mathcal{B}$ with coefficients in $\ZZ_{\geq 0}[v^{\pm1}]$.
\end{enumerate}

Note that all elements in a strongly positive basis are universally positive. Indeed, consider multiplying a basis element $B\in \mathcal{B}$ by a cluster monomial $X_1^{(a_1,a_2)}$ with $a_1,a_2$  sufficiently large to clear the denominator in the initial cluster expansion of $B$.  By properties (3) and (4) of a strongly positive basis
the product $BX_1^{(a_1,a_2)}$ is in $\ZZ_{\geq 0}[v^{\pm1}][X_1,X_2]$ and thus we see that a strongly positive basis element has non-negative coefficients in its initial cluster expansion. It then follows that $\mathcal{B}$ is universally positive because it is independent of the choice of an initial cluster.

Kimura and Qin \cite{KQ} showed the existence of bases for acyclic skew-symmetric quantum cluster algebras (hence for $\cA_v(b,c)$ with $b=c$) which satisfy (1), (3), (4).  In the rank 2 case, their bases are exactly the triangular bases, hence satisfy (2) as well (see \cite[Remark 1.8]{triangular}). 

\begin{conjecture}\label{str_pos_pointed} Given any strongly positive basis  $\mathcal{B}$
there exists a unique basis element pointed at $(a_1,a_2)$ for each $(a_1,a_2)\in\mathbb{Z}^2$.
\end{conjecture}

\noindent It is easy to see that property (3) of strongly positive bases together with Conjecture~\ref{str_pos_pointed} imply that every element of $\mathcal{B}$ is pointed.  We note that the conclusion of Conjecture~\ref{str_pos_pointed} holds for both the quantum greedy basis and the triangular basis.

In what follows, we write $Y\leq Z$ for $Y, Z\in\mathcal{A}_v(b,c)$ if $Z-Y$ is a Laurent polynomial with coefficients in $\mathbb{Z}_{\ge0}[v^{\pm1}]$.


\begin{conjecture}\label{upper_triangular}\mbox{}

{\rm(a)} There exists a unique strongly positive basis $\mathcal{B}^{upper}$ satisfying 
\begin{equation*}
\left(\begin{array}{c} \text{the basis element in}\\ \mathcal{B}  \text{ pointed at $(a_1,a_2)$} \end{array}  \right)
\leq \left(\begin{array}{c} \text{the basis element in}\\ \mathcal{B}^{upper}  \text{pointed at $(a_1,a_2)$}  \end{array}  \right)
\end{equation*}
for every strongly positive basis  $\mathcal{B}$ and every $(a_1,a_2)\in\mathbb{Z}^2$.

{\rm(b)} The triangular basis is $\mathcal{B}^{upper}$.
 \end{conjecture}

 \begin{conjecture}\label{conj:lower-greedy} Suppose that $b|c$ or $c|b$.

{\rm(a)} There exists a unique strongly positive basis $\mathcal{B}_{lower}$ satisfying \begin{equation*}
\left(\begin{array}{c} \text{the basis element in}\\ \mathcal{B}_{lower}  \text{ pointed at $(a_1,a_2)$} \end{array}  \right)
\leq \left(\begin{array}{c} \text{the basis element in}\\ \mathcal{B}  \text{ pointed at $(a_1,a_2)$} \end{array}  \right)
\end{equation*}
for every strongly positive basis  $\mathcal{B}$ and every $(a_1,a_2)\in\mathbb{Z}^2$.

{\rm(b)} The quantum greedy basis  is $\mathcal{B}_{lower}$.
 \end{conjecture}

In contrast to the example at the end of Section 3, we do not observe a failure of positivity  when $(b,c)=(1,5)$, $(1,6)$, $(2,4)$, $(2,6)$, $(3,3)$, $(3,6)$. This provides the motivation for assuming the condition `$b|c$ or $c|b$' in Conjecture~\ref{conj:lower-greedy}.

Conjectures~\ref{str_pos_pointed},~\ref{upper_triangular} and \ref{conj:lower-greedy} are trivially true for for finite types ($bc<4$), and are also known for affine types ($bc=4$) \cite{CDS, DX}.

Next we study the expansion coefficients relating the various bases. Let us begin by introducing notation.  For $(a_1,a_2),(a'_1,a'_2)\in\ZZ^2$ we define expansion coefficients $q^{a_1,a_2}_{a'_1,a'_2},r^{a_1,a_2}_{a'_1,a'_2}\in\ZZ[v^{\pm1}]$ as follows:
 \begin{align}
  \label{eq:greedy-standard expansion}
  X[a_1,a_2]&=M[a_1,a_2]+\sum\limits_{(a'_1,a'_2)<(a_1,a_2)} q^{a_1,a_2}_{a'_1,a'_2}M[a'_1,a'_2];\\
  \label{eq:triangular-greedy expansion}
  C[a_1,a_2]&=X[a_1,a_2]+\sum\limits_{(a'_1,a'_2)<(a_1,a_2)} r^{a_1,a_2}_{a'_1,a'_2}X[a'_1,a'_2],
 \end{align}
 where we write $(a'_1,a'_2)<(a_1,a_2)$ when $a'_1<a_1$ and $a'_2<a_2$.
 
  We now derive a recursion 
  for the expansion coefficients in equation \eqref{eq:triangular-greedy expansion} relating the triangular basis to the quantum greedy basis.  Our main ingredients will be the defining properties (P1) and (P2) of the triangular basis. Note that since each $X[a_1,a_2]$ is bar-invariant, we have $\overline{r^{a_1,a_2}_{a'_1,a'_2}}=r^{a_1,a_2}_{a'_1,a'_2}$ by (P1) .
 \medskip

 For $f\in\ZZ[v^{\pm1}]$ let $[f]_{\le0}$ denote the non-positive part of $f$, i.e. in $f$ we evaluate $v=0$ in all terms for which this makes sense.
 With this notation, we note that since $r^{a_1,a_2}_{a'_1,a'_2}$ is bar-invariant, it is determined by $[r^{a_1,a_2}_{a'_1,a'_2}]_{\le0}$.

 \begin{proposition}\label{rq_recursion}
    \begin{equation}\label{eq:recursion}
   \Big[r^{a_1,a_2}_{a''_1,a''_2}\Big]_{\le0}=-\sum\limits_{(a_1',a_2')>(a_1'',a_2'')} \Big[r^{a_1,a_2}_{a'_1,a'_2}q_{a_1'',a_2''}^{a_1',a_2'}\Big]_{\le0}.
  \end{equation}
\end{proposition}
\noindent This says that we have an easy recursive way to compute the decomposition in equation \eqref{eq:triangular-greedy expansion} if we know the decomposition in equation \eqref{eq:greedy-standard expansion}.

This proposition can be proved by reducing equation \eqref{eq:greedy-standard expansion} mod $\cL$ to get
 \begin{equation}\label{eq:f_red}
 \aligned
&  r^{a_1,a_2}_{a'_1,a'_2}X[a_1',a_2']\equiv \Big[r^{a_1,a_2}_{a'_1,a'_2}\Big]_{\le0}M[a_1',a_2']+\\
&\quad\sum\limits_{(a_1'',a_2'')<(a_1',a_2')} \Big[r^{a_1,a_2}_{a'_1,a'_2}q_{a_1'',a_2''}^{a_1',a_2'}\Big]_{\le0}M[a_1'',a_2'']\mod\cL.
\endaligned
 \end{equation}
 Keeping in mind (P2), we reduce the equality \eqref{eq:triangular-greedy expansion} mod $\cL$ and apply equation \eqref{eq:f_red} to arrive at the desired recursion.
\medskip

Based on extensive computations using Proposition~\ref{rq_recursion} we make the following positivity conjecture.

 \begin{conjecture}\label{CXcoeff:pos}
  For any $(a_1,a_2)\in\ZZ^2$ the expansion coefficients of the triangular basis element $C[a_1,a_2]$ in terms of the quantum greedy basis are positive; more precisely, we have

 $$r^{a_1,a_2}_{a'_1,a'_2}\in\ZZ_{\ge0}[v^{\pm1}]\quad\text{for all}\quad (a'_1,a'_2)\in\ZZ^2.$$

More generally,  the expansion coefficients of any strongly positive basis in terms of the quantum greedy basis are positive.
 \end{conjecture}

We end this article by summarizing the aforementioned proven or conjectured properties of   the quantum greedy basis and the triangular basis.

\begin{tabular}{|l|c|}
\hline
&quantum greedy basis\\
\hline
bar-invariant& by Proposition~\ref{recursive-definition-q}\\
independent of initial cluster&new result\\
contains all cluster monomials&new result\\
positive expansion of products& conjectured \\
consists of indecomposables&true if universally positive \\
each element is pointed & by definition \\
support for $(a_1,a_2)\in\Phi^{im}_+$ & Figure 1\\
bound of strongly positive bases&conjectured\\
\hline
\end{tabular}
\begin{center}
Table 1. Properties of the quantum greedy basis for $b|c$ or $c|b$
\end{center}

\begin{tabular}{|l|c|}
\hline
&triangular basis\\
\hline
bar-invariant& by definition\\
independent of initial cluster&proved in \cite{triangular}\\
contains all cluster monomials&proved in \cite{triangular}\\
positive expansion of products& conjectured \\
& (proved for $b=c$ in \cite{KQ} )\\
consists of indecomposables&false for $bc\geq4$  \\
each element is pointed & by definition \\
support for $(a_1,a_2)\in\Phi^{im}_+$& conjectured, see Figure 2\\
bound of strongly positive bases&conjectured\\
\hline
\end{tabular}
\begin{center}
Table 2. Properties of the triangular basis
\end{center}

\begin{acknowledgments}
Most of the ideas toward this work from D.~R. and A.~Z. were had during their stay at the Mathematical Sciences Research Institute as part of the Cluster Algebras Program.  They would like to thank the MSRI for their hospitality and support.
Research supported in part by NSF grants DMS-0901367 (K.~L.) and DMS-1103813 (A.~Z.),  and by Oakland University URC Faculty Research Fellowship Award (L.~L.).
The authors would like to thank F.~Qin for valuable dicussions. They would also like to thank S.~Fomin for many constructive suggestions.
\end{acknowledgments}





\end{article}








\end{document}